\documentclass[11pt,a4paper]{article}

\usepackage[hmargin=24mm,vmargin=22mm]{geometry}
\usepackage{newtxtext}
\usepackage{amsmath,amsthm}
\usepackage{newtxmath}
\usepackage{microtype}

\usepackage{xcolor}
\usepackage[unicode,colorlinks=true,allcolors=blue!80!black]{hyperref}

\newtheorem{theorem}{Theorem}
\newtheorem{lemma}{Lemma}
\newtheorem{corollary}{Corollary}
\DeclareMathOperator{\FC}{FC}

\title{Exact factorizations of countable groups\\
into discrete and thin subsets}
\author{Mikhail Kabenyuk}
\date{}

\begin{document}
\maketitle

\begin{abstract}
For every countably infinite Hausdorff topological group $G$, we
construct an exact factorization $G=AB$, where $A,B\subseteq G$
and every element of $G$ has a unique representation as $ab$
with $a\in A$ and $b\in B$. Both factors are infinite, closed,
discrete, and thin. A subset $C\subseteq G$ is thin if the
intersections $hC\cap C$ and $Ch\cap C$ are finite for every
nonidentity element $h\in G$.

In any countably infinite group, an exact factorization into two
infinite thin subsets can be constructed starting from any partial
exact factorization with finite nonempty factors. Given a finite or
countably infinite family of infinite subsets of the group, we can
choose the extension so that each factor meets every member of the
family in infinitely many points. If the topology on the group has
a countable base and every nonempty open set is infinite, taking
the nonempty members of the base as this family yields dense thin
factors.
\end{abstract}

\section{Introduction}

Let $G$ be a group with identity $e$. The product $AB$ of subsets
$A,B\subseteq G$ is called a \textit{partial exact factorization}
if the map $(a,b)\mapsto ab$ is injective on $A\times B$. If, in
addition, $AB=G$, we call it an \textit{exact factorization} of $G$.

Protasov \cite{Protasov2016} proved that every countable Hausdorff
topological group $G$ has a closed discrete subset $C$ such that
$G=CC^{-1}$. He also observed that $C$ can be chosen to be left
thin, meaning that $hC\cap C$ is finite for every $h\ne e$.
In the same paper, he asked whether such a group admits an exact
factorization into two discrete subsets and whether both subsets
can be chosen to be closed. We answer both questions affirmatively.

\begin{theorem}\label{thm:discrete}
Every countably infinite Hausdorff topological group $G$ admits an
exact factorization $G=AB$ into two infinite closed discrete subsets
$A$ and $B$.
\end{theorem}

Following Protasov \cite{Protasov2013}, we call a subset
$C\subseteq G$ \textit{thin} if
\begin{equation}\label{eq:thin}
  |hC\cap C|<\infty,\quad |Ch\cap C|<\infty
  \quad (h\in G\setminus\{e\}).
\end{equation}

Lutsenko \cite[Theorem~3]{Lutsenko2010} proved that every infinite
group $G$ has a subset $C$ such that $G=CC^{-1}$ and each of the
intersections in \eqref{eq:thin} contains at most four elements.
Protasov \cite{Protasov2021} obtained an exact factorization into
two infinite thin subsets for every countably infinite group
that admits a nondiscrete Hausdorff group topology. 

Protasov \cite{Protasov2013} proved that every countably infinite
precompact Hausdorff topological group has a dense thin subset.
Protasov and Slobodianiuk \cite{ProtasovSlobodianiuk2017} studied
exact factorizations into dense subsets. They proved that a group
$G$ of cardinality $\kappa>\aleph_0$ admits an exact factorization
into two dense subsets of cardinality $\kappa$ if the topology on
$G$ has a base of cardinality at most $\kappa$ and every nonempty
open set has cardinality $\kappa$.

The countable case follows from the extension theorem
\cite[Theorem~1.1 and Corollary~5.2]{Kabenyuk2026}. Given an at most
countable family of infinite subsets of the group, this theorem
extends any partial exact factorization with finite nonempty factors
to an exact factorization in which each factor meets every member
of the family in infinitely many points. Theorem~\ref{thm:thin}
shows that both factors can also be chosen to be thin and extends
Protasov's result in \cite{Protasov2021} to all countably infinite
groups.

\begin{theorem}\label{thm:thin}
Let $G$ be a countably infinite group and let $\mathcal X$ be an
at most countable family of infinite subsets of $G$. Every partial
exact factorization $A_0B_0$ with finite nonempty factors extends
to an exact factorization $G=AB$ into two infinite thin subsets
$A$ and $B$ satisfying
\[
  A_0\subseteq A,\quad B_0\subseteq B,\quad
  |A\cap X|=|B\cap X|=\aleph_0 \quad (X\in\mathcal X).
\]
If $A_0$ and $B_0$ are singletons, $A$ and $B$ can be chosen so that
for each of them both intersections in \eqref{eq:thin} contain at
most two elements.
\end{theorem}

If the topology on $G$ has a countable base and every nonempty open
set is infinite, Theorem~\ref{thm:thin} gives an exact factorization
into dense thin subsets (Corollary~\ref{cor:dense-thin}).

We prove Theorem~\ref{thm:discrete} in Section~2. As in
\cite{Protasov2016}, we treat precompact and nonprecompact groups
separately. In Section~3, we prove Theorem~\ref{thm:thin} and show
that the closed discrete factors in Theorem~\ref{thm:discrete} can
also be chosen to be thin. To ensure thinness in both constructions,
we use a variant of a technique of Chou \cite{Chou1969}. When adding
a point $z$ to the current factor $C$, we choose it outside
$CC^{-1}C$.

\section{Closed discrete factors}

For a group $L$ and an element $h\in L$, put
$R_L(h)=\{z\in L:z^2=h\}$. The extension step in the proof of
\cite[Lemma~5.1]{Kabenyuk2026} amounts to excluding the elements
described in the following lemma.

\begin{lemma}\label{lem:extension}
Let $AB$ be a partial exact factorization with finite nonempty
factors and let $g\in G\setminus AB$. Put
\begin{align*}
  E&=ABB^{-1}\cup gB^{-1}A^{-1}A,\\
  D&=\bigcup_{(a,b)\in A\times B} R_G(abg^{-1})\,gb^{-1}.
\end{align*}
If $x\in G\setminus(E\cup D)$ and $y=x^{-1}g$, then
$x\notin A$, $y\notin B$, and $(A\cup\{x\})(B\cup\{y\})$
is a partial exact factorization containing $g$. 
\end{lemma}

\begin{proof}
Since $A\subseteq ABB^{-1}$ and
$gB^{-1}\subseteq gB^{-1}A^{-1}A$, we have $x\notin A$
and $y\notin B$. The set of all products is the union of
\[
  AB,\quad xB,\quad Ay,\quad \{g\}.
\]
Within each of these four sets, the corresponding representation
is unique. The equalities $xb=a'b'$ and $ay=a'b'$ imply
$x=a'b'b^{-1}$ and $x=g(b')^{-1}(a')^{-1}a$, respectively,
and both are ruled out by $x\notin E$.
By assumption, $g\notin AB$, while $g=xb$ or $g=ay$ would give
$y=b$ or $x=a$, respectively. Finally,
\[
  xb=ay
  \quad\Longleftrightarrow\quad xbg^{-1}x=a
  \quad\Longleftrightarrow\quad (xbg^{-1})^2=abg^{-1}.
\]
Thus $x\notin D$ ensures that $xB$ and $Ay$ are disjoint.
\end{proof}

\noindent\textbf{Remark.} Note that the elements $abg^{-1}$ in the
definition of $D$ are pairwise distinct by exactness of $AB$, and
$g\notin AB$ implies that none of them equals $e$.

\begin{lemma}\label{lem:haar}
Let $K$ be an infinite compact Hausdorff group with normalized
Haar measure $\mu$, and let $G\le K$ be a countable subgroup. Then
\[
  \mu(S)<1,\quad S=\{z\in K:z^2\in G\setminus\{e\}\}.
\]
\end{lemma}

\begin{proof}
The set $S$ is the countable union of the closed sets $R_K(h)$,
$h\in G\setminus\{e\}$, and is therefore Borel. Suppose that
$\mu(S)=1$. Then $R=R_K(h)$ has positive measure for some
$h\in G\setminus\{e\}$.

For each $t\in K$, the set $S\cap t^{-1}S$ has full measure,
and every $x$ in this set satisfies
\[
  xtx^{-1}=t^{-1}(tx)^2(x^2)^{-1}\in t^{-1}G.
\]
Since $t^{-1}G$ is countable and $xtx^{-1}=yty^{-1}$ is equivalent
to $x\in yC_K(t)$, the set $S\cap t^{-1}S$ is covered by countably
many left cosets of the centralizer $C_K(t)$. Hence
$\mu(C_K(t))>0$. 
Thus $C_K(t)$ has finite index in $K$ and, being closed, is open.

The open sets $C_K(z)$, $z\in R$, cover the compact set $R$, since
$z\in C_K(z)$. Taking a finite subcover, we find that
$\mu(R\cap C_K(z))>0$ for some $z\in R$. For $x\in R\cap C_K(z)$
we have $(z^{-1}x)^2=z^{-2}x^2=e$. Thus $z^{-1}(R\cap C_K(z))$
has positive measure and is disjoint from $S$, contradicting
$\mu(S)=1$.
\end{proof}

\begin{lemma}\label{lem:nonprecompact-cosets}
Let $G$ be a Hausdorff topological group.
\begin{enumerate}
\renewcommand{\theenumi}{\textup{(\roman{enumi})}}
\renewcommand{\labelenumi}{\theenumi}
\item\label{lem:nonprecompact-translates}
If $G$ is nonprecompact, there is a symmetric open neighborhood
$U$ of the identity with infinitely many pairwise disjoint left
translates.

\item\label{lem:nonprecompact-cover}
Let $U$ be any symmetric open neighborhood of the identity with
infinitely many pairwise disjoint left translates. Then for every
finite $F\subseteq G$ and every finite family of left cosets
$c_kH_k$ of subgroups of infinite index, we have
\[
  FU\cup\bigcup_{k=1}^q c_kH_k\ne G.
\]
\end{enumerate}
\end{lemma}

\begin{proof}
Part~\ref{lem:nonprecompact-translates}.
Let $G$ be a nonprecompact Hausdorff topological group.
By the criterion for precompactness (see, for example,
\cite[Section~3.7, Corollary~3.7.17]{ArhangelskiiTkachenko2008}),
there is a symmetric open neighborhood $W$ of the identity
such that $FW\ne G$ for every finite $F\subseteq G$.
Choose a symmetric open neighborhood $U$ of the identity with
$U^2\subseteq W$.

Set $p_1=e$ and choose $p_{n+1}\notin\bigcup_{i=1}^n p_iW$
recursively. By symmetry of $W$, we have $p_i^{-1}p_j\notin W$
for $i\ne j$. Since $UU^{-1}\subseteq W$, the sets $p_iU$
are pairwise disjoint.

For part~\ref{lem:nonprecompact-cover}, take any $U$ with the stated
properties and choose pairwise disjoint left translates
$p_1U,p_2U,\ldots$.
Fix $F$ and the cosets $c_kH_k$, and put $N=|F|+1$.
By Neumann's lemma \cite[Lemma~4.1]{Neumann1954}, we can choose
$s\in G$ outside all the sets $c_kH_kp_i^{-1}$ with
$1\le i\le N$ and $1\le k\le q$, since each of these sets is a
left coset of the infinite-index subgroup $p_iH_kp_i^{-1}$.
Then the points $sp_1,\ldots,sp_N$ lie outside all the cosets $c_kH_k$. 
Each set $fU$, $f\in F$, contains at most one of these points.
Indeed, suppose that $sp_i,sp_j\in fU$ with $i\ne j$.
Then $p_i^{-1}p_j\in U^{-1}U=UU^{-1}$ by symmetry of $U$.
Writing $p_i^{-1}p_j=uv^{-1}$ with $u,v\in U$, we obtain
$p_i u=p_j v\in p_iU\cap p_jU$, a contradiction.
Hence at least one of the $N$ points also lies outside $FU$.
\end{proof}

\begin{lemma}\label{lem:nonprecompact}
Let $G$ be a nonprecompact Hausdorff topological group. There is
a symmetric open neighborhood $U$ of the identity with the
following property. If $F\subseteq G$ is finite,
$h_1,\ldots,h_r\in G\setminus\{e\}$ are pairwise distinct,
$r\ge 1$, and $t_1,\ldots,t_r\in G$, then
\[
  FU\cup \bigcup_{i=1}^r R_G(h_i)t_i \ne G.
\]
\end{lemma}

\begin{proof}
Take $U$ as in
Lemma~\ref{lem:nonprecompact-cosets}\ref{lem:nonprecompact-translates}.
Suppose that
$G=FU\cup D$, where $D=\bigcup_{i=1}^r R_G(h_i)t_i$.
Put $u_i=t_i^{-1}$ and $v_i=h_it_i$. Then $x\in R_G(h_i)t_i$
if and only if $xu_ix=v_i$.

For every $s\in G$ and every $x\notin(F\cup s^{-1}F)U$, both $x$ and $sx$ belong to $D$. 
Hence for some $i,j$ we have
\[
  xu_ix=v_i,\quad (sx)u_j(sx)=v_j.
\]
Substituting $u_i^{-1}x^{-1}v_i$ for the final factor $x$ in
$xu_jsx=s^{-1}v_j$ gives
\[
  x(u_jsu_i^{-1})x^{-1}=s^{-1}v_jv_i^{-1}.
\]
Thus $G$ is covered by $(F\cup s^{-1}F)U$ together with the
solution sets of these equations as $i$ and $j$ vary. 
Each nonempty solution set is a left coset of the centralizer
$C_G(u_jsu_i^{-1})$. 
By Lemma~\ref{lem:nonprecompact-cosets}\ref{lem:nonprecompact-cover},
at least one of these cosets is a coset of a centralizer of
finite index.

It follows that for each $s\in G$ there are $i,j$ such that
$s\in u_j^{-1}\FC(G)u_i$, where $\FC(G)$ is the normal subgroup
of elements with finite conjugacy classes. The sets
$u_j^{-1}\FC(G)u_i$ are cosets of $\FC(G)$, so
$[G:\FC(G)]\le r^2$. By \cite[Lemma~4.1]{Kabenyuk2026}, 
the group $G$ is amenable as an abstract group.

Let $\mu$ be a bi-invariant finitely additive probability measure
on all subsets of $G$. By \cite[Proposition~3.3]{Kabenyuk2026},
we have $\mu(D)<1$. 
By Lemma~\ref{lem:nonprecompact-cosets}\ref{lem:nonprecompact-translates},
for every $n\ge 1$ there are $n$ pairwise disjoint left translates
of $U$, so $n\mu(U)\le 1$.
Thus $\mu(U)=0$ and $\mu(FU)=0$.
This gives the contradiction
\[
  1=\mu(G)\le \mu(FU)+\mu(D)<1.\qedhere
\]
\end{proof}

\begin{proof}[Proof of Theorem~\textup{\ref{thm:discrete}}]
Enumerate $G$ as $g_0,g_1,\ldots$ with $g_0=e$, and put
$A_0=B_0=\{e\}$.
Suppose that $A_n$ and $B_n$ have been constructed and that
$A_nB_n$ is a partial exact factorization. If $g_n\in A_nB_n$,
put $A_{n+1}=A_n$ and $B_{n+1}=B_n$.
If $g_n\notin A_nB_n$, choose $x\in G\setminus(E\cup D)$, where
$E$ and $D$ are as in Lemma~\ref{lem:extension} for $A_n,B_n,g_n$,
and put
\[
  A_{n+1}=A_n\cup\{x\},\quad B_{n+1}=B_n\cup\{x^{-1}g_n\}.
\]
Below we treat the precompact and nonprecompact cases separately,
specifying the additional conditions on $x$ and showing that they
can be satisfied.

At every stage, $A_{n+1}B_{n+1}$ is a partial exact factorization
containing $g_n$.
Set $A=\bigcup_{n\ge 0}A_n$ and $B=\bigcup_{n\ge 0}B_n$.

\medskip
\noindent\textit{Precompact case.}
Suppose that $G$ is precompact. 
Identify it with a dense subgroup of its compact completion $K$. 
Let $\mu$ be normalized Haar measure on $K$ and put $\delta=1-\mu(S)>0$, 
where $S$ is the set in Lemma~\ref{lem:haar}.

At each stage $n\ge 0$, once $A_n$ and $B_n$ have been constructed,
choose a closed neighborhood $U_n$ of $g_n$ in $K$ such that
\[
  U_n\cap(A_n\cup B_n)\subseteq\{g_n\},\quad
  \mu(U_n)<\frac{\delta}{2^{n+3}}.
\]
To see that such a $U_n$ exists, note that every singleton in an
infinite compact group has measure zero.
By outer regularity of Haar measure, $g_n$ has an open neighborhood
of measure less than $\delta/2^{n+3}$. 
Remove $(A_n\cup B_n)\setminus\{g_n\}$ from this neighborhood. 
By regularity of $K$ as a topological space, 
we can choose a closed neighborhood $U_n$ contained in the resulting open set.

When adding a new pair, we require both $x$ and $x^{-1}g_n$ to lie
outside $P_n=\bigcup_{i=0}^n U_i$. 
To ensure this, it suffices to exclude $P_n\cup g_nP_n^{-1}$. 
This set is closed since $P_n$ is closed and inversion and left translation are homeomorphisms.
By invariance of normalized Haar measure under translations
and inversion (see, for example, \cite[Theorem~9.3.13,
items~5--7 and p.~590]{ArhangelskiiTkachenko2008}), we have
\[
  \mu(P_n\cup g_nP_n^{-1})
  \le \mu(P_n)+\mu(g_nP_n^{-1})
  =2\mu(P_n)<\frac{\delta}{2}.
\]

If $g_n\notin A_nB_n$, also consider the subset of $K$ given by
\[
  D_K=\bigcup_{(a,b)\in A_n\times B_n} R_K(abg_n^{-1})\,g_n b^{-1}.
\]
It is closed and satisfies $D_K\cap G=D$ for the set $D$ in Lemma~\ref{lem:extension}. 
The elements $abg_n^{-1}$ are pairwise distinct and lie in $G\setminus\{e\}$.
By right invariance of Haar measure,
\[
  \mu(D_K)
  \le \sum_{(a,b)\in A_n\times B_n}\mu\bigl(R_K(abg_n^{-1})\bigr)
  \le \mu(S).
\]
The last inequality holds because the sets $R_K(abg_n^{-1})$
are pairwise disjoint subsets of $S$.

The set $E$ is finite, so it is closed in $K$ and $\mu(E)=0$.
The union $E\cup D_K\cup P_n\cup g_nP_n^{-1}$ is closed and has
measure less than $\mu(S)+\delta/2=1-\delta/2$.
Its open complement is nonempty and therefore meets $G$.
Choose $x$ in this intersection.

For every $n$ we have $U_n\cap(A\cup B)\subseteq\{g_n\}$.
Thus each point of $G$ has a neighborhood whose only possible
point in either factor is the point itself. 
It follows that $A$ and $B$ are closed and discrete in $G$.

\medskip
\noindent\textit{Nonprecompact case.}
Now suppose that $G$ is nonprecompact. Take $U$ from
Lemma~\ref{lem:nonprecompact}. When adding a new pair, we require
\[
  x\notin A_nU,\quad x^{-1}g_n\notin UB_n.
\]
By symmetry of $U$, the second condition is equivalent to
$x\notin g_nB_n^{-1}U$. Since $A_n$ and $g_nB_n^{-1}$ are both
contained in $E$, it suffices to choose $x$ outside $EU\cup D$,
with $E$ and $D$ as in Lemma~\ref{lem:extension}. This is possible
by Lemma~\ref{lem:nonprecompact}.

It follows that $a^{-1}a'\notin U$ for distinct $a,a'\in A$
and $b'b^{-1}\notin U$ for distinct $b,b'\in B$.
Choose a symmetric open neighborhood $V$ of the identity with
$V^2\subseteq U$. For any $z\in G$, the neighborhood $zV$ contains
at most one point of $A$, and $Vz$ contains at most one point of
$B$. If $z$ lies outside the relevant factor, deleting any such
point leaves an open neighborhood of $z$. Thus both factors are
closed and discrete.

\medskip
In both the precompact and nonprecompact cases, the construction
gives $g_n\in A_{n+1}B_{n+1}$ for every $n$, so $G=AB$.
If an element had two distinct
representations, both would already occur at some finite stage,
contrary to the exactness of $A_nB_n$ at that stage. Since $G=AB$
is infinite, new pairs were added at infinitely many stages.
Both factors grew at every such stage, so $A$ and $B$ are infinite.
\end{proof}

\section{Thin factors}

\begin{lemma}\label{lem:thin-extension}
Let $C$ be a finite nonempty subset of a group $G$, let
$z\in G\setminus CC^{-1}C$, and put $C'=C\cup\{z\}$. 
Then for every $h\in G\setminus\{e\}$ we have
\[
\begin{aligned}
  hC\cap C\ne\varnothing &\implies hC'\cap C'=hC\cap C,\\
  hC\cap C=\varnothing   &\implies hC'\cap C'\subseteq\{z,hz\},\\[2pt]
  Ch\cap C\ne\varnothing &\implies C'h\cap C'=Ch\cap C,\\
  Ch\cap C=\varnothing   &\implies C'h\cap C'\subseteq\{z,zh\}.
\end{aligned}
\]
\end{lemma}

\begin{proof}
If $hC\cap C\ne\varnothing$, then $h,h^{-1}\in CC^{-1}$.
In any new pair of the form $c,hc$, one point must be $z$ and
the other must lie in $C$, since $h\ne e$. This would give
$z\in hC\cup h^{-1}C\subseteq CC^{-1}C$, a contradiction.
Similarly, if $Ch\cap C\ne\varnothing$, then $h,h^{-1}\in C^{-1}C$,
and a new pair of the form $c,ch$ would give
$z\in Ch\cup Ch^{-1}\subseteq CC^{-1}C$. If the original
intersection was empty, then after adjoining $z$ it can contain
only $z$ and $hz$ for the left translate, or only $z$ and $zh$
for the right translate.
\end{proof}

\begin{proof}[Proof of Theorem~\textup{\ref{thm:thin}}]
By \cite[Lemma~5.1]{Kabenyuk2026}, in any infinite group, given
a partial exact factorization $AB$ with finite nonempty factors
and an element $g\notin AB$, one can add $x$ to $A$ and $x^{-1}g$
to $B$ while preserving exactness. Moreover, $x$ can be chosen
outside any prescribed finite set.

Enumerate $G=\{g_0,g_1,\ldots\}$. Choose a sequence $(X_n)_{n\ge 0}$
of members of $\mathcal X\cup\{G\}$ in which every member occurs
infinitely often. Starting with $A_0,B_0$, construct increasing
sequences of finite sets $A_n,B_n$ so that $A_nB_n$ remains a
partial exact factorization at each stage.

At stage $n$, choose $a_n$ and then $b_n$ as follows:
\begin{align*}
  a_n &\in X_n\setminus
    \bigl(A_nB_nB_n^{-1}\cup A_nA_n^{-1}A_n\bigr),
  & \widetilde A_n &= A_n\cup\{a_n\},\\
  b_n &\in X_n\setminus
    \bigl(\widetilde A_n^{-1}\widetilde A_n B_n\cup B_nB_n^{-1}B_n\bigr),
  & \widetilde B_n &= B_n\cup\{b_n\}.
\end{align*}
These choices are possible since $X_n$ is infinite and all the
excluded sets are finite. The inclusion $C\subseteq CC^{-1}C$
gives $a_n\notin A_n$ and $b_n\notin B_n$. Excluding
$A_nB_nB_n^{-1}$ and $\widetilde A_n^{-1}\widetilde A_n B_n$
ensures that $a_n B_n$ is disjoint from $A_n B_n$ and that
$\widetilde A_n b_n$ is disjoint from $\widetilde A_n B_n$,
respectively. Hence $\widetilde A_n\widetilde B_n$ is a partial
exact factorization.

If $g_n\in\widetilde A_n\widetilde B_n$, put
$A_{n+1}=\widetilde A_n$ and $B_{n+1}=\widetilde B_n$.
Otherwise, apply \cite[Lemma~5.1]{Kabenyuk2026} to choose $x$
subject to the additional condition
\[
  x\notin \widetilde A_n\widetilde A_n^{-1}\widetilde A_n
  \,\cup\, g_n\bigl(\widetilde B_n\widetilde B_n^{-1}\widetilde B_n\bigr)^{-1}.
\]
Then $y=x^{-1}g_n\notin\widetilde B_n\widetilde B_n^{-1}\widetilde B_n$.
Put $A_{n+1}=\widetilde A_n\cup\{x\}$ and
$B_{n+1}=\widetilde B_n\cup\{y\}$.

Let $A=\bigcup_{n\ge 0}A_n$ and $B=\bigcup_{n\ge 0}B_n$.
Each $g_n$ belongs to $A_{n+1}B_{n+1}$, and any two distinct
representations of the same element would already occur at a
finite stage. Hence $G=AB$ is an exact factorization.
For every $X\in\mathcal X\cup\{G\}$, there are infinitely many
stages with $X_n=X$. At each such stage, a new point of $X$ is added
to each factor. Consequently, $A\cap X$ and $B\cap X$ are infinite.
Taking $X=G$ shows that $A$ and $B$ are infinite.

Consider one of the resulting factors $C\in\{A,B\}$ and let $D$
denote its current finite approximation. Each new point is chosen
outside $DD^{-1}D$ (cf.~\cite[proof of Proposition~4.1,
p.~202]{Chou1969}). By Lemma~\ref{lem:thin-extension}, for a fixed
$h\ne e$, each of the intersections $hD\cap D$ and $Dh\cap D$
remains unchanged once it first becomes nonempty. Since any two
points of $C$ belong to some finite approximation, $hC\cap C$
and $Ch\cap C$ are finite. Thus $C$ is thin. If $A_0$ and $B_0$
are singletons, both intersections are empty in the initial
approximation. By the same lemma, each contains at most two points
when it first becomes nonempty. Therefore $|hC\cap C|\le 2$
and $|Ch\cap C|\le 2$.
\end{proof}

\begin{corollary}\label{cor:dense-thin}
Suppose that the topology on a countably infinite group $G$
has a countable base and that every nonempty open set is infinite.
Every partial exact factorization with finite nonempty factors
extends to an exact factorization $G=AB$ into dense thin subsets.
Each factor meets every nonempty open set in infinitely many
points.
\end{corollary}

\begin{proof}
In Theorem~\ref{thm:thin}, take $\mathcal X$ to be the family
of nonempty members of a countable base. Each nonempty open set
contains one of these members.
\end{proof}

\begin{corollary}
The factors in Theorem~\ref{thm:discrete} can be chosen so that
each is thin and each of its intersections in \eqref{eq:thin}
has at most two elements.
\end{corollary}

\begin{proof}
When adding a pair in the proof of Theorem~\ref{thm:discrete},
we also require $x$ to lie outside the finite set
\[
  Q_n=A_nA_n^{-1}A_n\cup g_n(B_nB_n^{-1}B_n)^{-1}.
\]
In the precompact case, $Q_n$ is closed in $K$ and has measure
zero. In the nonprecompact case, apply Lemma~\ref{lem:nonprecompact}
with $F=E\cup Q_n$. Then $Q_n\subseteq FU$ since $e\in U$.
In both cases, every point added to the current factor $C$ is
chosen outside $CC^{-1}C$. Since $A_0=B_0=\{e\}$,
Lemma~\ref{lem:thin-extension} gives the bound of two for both
intersections, as at the end of the proof of Theorem~\ref{thm:thin}.
\end{proof}

\end{document}